\documentclass{article}
\usepackage{maa-monthly}
\usepackage{tikz}
\usepackage{float}
\definecolor[named]{Shade}{RGB}{215,215,215}
\definecolor[named]{Lightsilver}{RGB}{240,240,240}
\definecolor[named]{MyGray}{RGB}{170,170,170}

\usetikzlibrary{fadings}
\pgfdeclareradialshading{myring}{\pgfpointorigin}
{
color(0cm)=(transparent!50);
color(4mm)=(pgftransparent!75);
color(8mm)=(pgftransparent!100)
}
\pgfdeclarefading{ringo}{\pgfuseshading{myring}}

\theoremstyle{theorem}
\newtheorem{theorem}{Theorem}
\newtheorem{lemma}{Lemma}
 
\theoremstyle{definition}

\begin{document}

\title{A New View of Hypercube Genus}
\markright{Hypercube Genus}
\author{Richard H.~Hammack and Paul C.~Kainen}

\maketitle

\begin{abstract}
Beineke, Harary and Ringel discovered a formula for the minimum genus of a
torus in which the $n$-dimensional hypercube graph can be embedded.  
We give a new proof of the formula by building this surface
as a union of certain faces in the hypercube's 2-skeleton. For odd dimension $n$, the entire 2-skeleton decomposes into  
$(n-1)/2$ copies of the surface, and the intersection of any two copies is
the hypercube graph.
\end{abstract}

\section{introduction.}
\begin{figure}[b]
\centering
\includegraphics{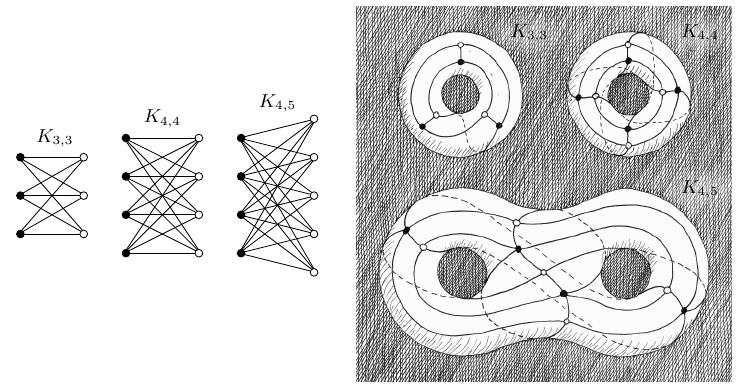}
\caption{ Left: The complete bipartite graph $K_{m,n}$ can be regarded as having $m$ black vertices,
$n$ white vertices, and an edge joining any two vertices of different colors.
Right: $K_{3,3}$,  $K_{4,4}$, and  $K_{4,5}$ on surfaces.}
\label{Fig:K55}
\end{figure}

\noindent
What graphs can be drawn on what surfaces without crossed edges?
Kuratowski's theorem~\cite[Theorem 6.18]{CHARTRAND} implies that
the complete bipartite graph $K_{m,n}$ (Figure~\ref{Fig:K55})
cannot
be drawn on the sphere (or plane) without crossed edges if $\min\{m,n\}\geq 3$.
One can try to draw $K_{3,3}$ on the sphere, but will always fail.
However, it can be drawn on the torus, as shown in Figure~\ref{Fig:K55}.

The {\bf genus} of a graph $G$, denoted $\gamma(G)$, is the least integer $g$
for which $G$ can be drawn on a closed, connected,  orientable surface of genus $g$ without edges crossing.
(A sphere has genus 0; a surface with $g$ holes has genus $g$.) Thus $\gamma(K_{3,3})=1$.

Figure~\ref{Fig:K55} also shows $K_{4,4}$ and $K_{4,5}$, together with
drawings of them on tori of genus~1 and 2, respectively.
Indeed, Ringel~\cite{RINGEL2} proved that for $2 \leq m \leq n$,
\[\gamma(K_{m,n})=\left\lceil \frac{(m-2)(n-2)}{4}\right\rceil.\]
Thus $\gamma(K_{4,4})=1$ and $\gamma(K_{4,5})=2$, so our drawings
in Figure~\ref{Fig:K55} are optimal.

Such genus formulas have been established for most well-known families of graphs. 
For instance, the difficult proof of the simple formula $\gamma(K_n) = \lceil (n-3)(n-4)/12\rceil$
was instrumental in settling the Heawood map coloring 
conjecture~\cite[Chapter 7]{CHARTRAND}, \cite{RINGEL3}.

For the $n$-dimensional hypercube graph $Q_n$, Ringel~\cite{RINGEL1} and Beineke and Harary~\cite{BEINEKE} used recursive arguments exploiting the hypercube's product structure
to deduce
\begin{equation}
\gamma(Q_n)=1+(n-4)2^{n-3}.
\label{Eqn:Qn}
\end{equation}
Their proofs lead to generalizations like 
\cite{HUNTER,pck-st-cr,k-white,PISANSKI1,PISANSKI2,WHITE}.  
In contrast, our short, visual proof directly constructs the genus surface from the square faces of
the $n$-cube and further shows that for odd $n$, the $2$-skeleton of the $n$-cube is the union of $(n{-}1)/2$ copies of the genus surface,  with no common faces, intersecting pairwise at $Q_n$.

Sections 2 and 3 define hypercubes, skeleta, and 2-cell embeddings of graphs in surfaces, and state Euler's
formula for genus.
Our proof of equation (1) is in Section 4
while Section 5 gives the bonus: a factorization of the $2$-skeleton into genus surfaces.

\section{Hypercubes.}
\label{Section:Hypercubes}

Let $I$ denote the unit interval $[0,1]$ and let $O$ denote its boundary $O := \partial I = \{0,1\}$.
(We use the notation $O$  because it will be
useful to regard an
interval as being ``active'' [$I$] or ``inactive'' [$O$] in the manner described below.)

The {\bf $n$-dimensional hypercube}, or {\bf $n$-cube} is the polytope
$H_n=I^n \subseteq\mathbb{R}^n$. Thus
$H_n=\{(x_1,\ldots, x_n) \mid 0\leq x_i\leq 1\}\subseteq \mathbb{R}^n$ is the 
intersection of the $2n$ half-spaces $x_i\geq 0$ and $x_i\leq1$, for $1\leq i\leq n$.
(See~\cite{ZIEGLER} for an introduction to polytopes.)

The $2^n$ {\bf vertices} of $H_n$ are the elements of $O^n$,
which we identify with the binary strings of length $n$. (For example,
$(1,0,1,0)$ is 1010, etc.)
The {\bf edges} of $H_n$ are the connected components of the products
$O\times O\times \cdots \times I\times \cdots\times O$
having one active factor $I$ and $n{-}1$ inactive factors $O$.
Thus $H_n$ has $n2^{n-1}$ edges, and each is a line segment joining two vertices that differ in
exactly one coordinate (namely the active coordinate).
The {\bf faces} of $H_n$ are the squares that are the connected components of 
\[O\times \cdots \times I\times  \cdots \times I\times \cdots\times O,\]
where two of the factors are $I$'s and the rest are $O$'s. Thus $H_n$ has $\binom{n}{2}2^{n-2}$ faces, and the boundary of each
face consists of four edges. Likewise $H_n$ has $\binom{n}{3}2^{n-3}$ {\bf 3-faces}
\[O\times \cdots \times I\times  \cdots \times I\times  \cdots \times I\times \cdots\times O,\]
formed by choosing three positions for the  $I$'s. Each 3-face is a 3-cube whose boundary has six faces.
In general, for each $0\leq k\leq n$, the $n$-cube has $\binom{n}{k}2^{n-k}$ {\bf $k$-faces} formed by
choosing $k$  active factors. Each $k$-face is a $k$-cube.
(For brevity we call the 0-faces, 1-faces, and 2-faces {\it vertices, edges, and faces}, respectively,
as stated above.)

The $(n{-}1)$-faces are called the {\bf facets} of $H_n$.
For each $1\leq k\leq n$, there are two {\it opposite} facets that are the
two components of
 $I\times I\times \cdots\times O\times \cdots \times I$, 
where the sole
$O=\{0,1\}$ is the $k$th factor. The reader should verify that any two nonopposite facets intersect in an
$(n{-}2)$-face. Collectively the facets form the boundary of $H_n$.

\begin{figure}[t]
\centering
\includegraphics{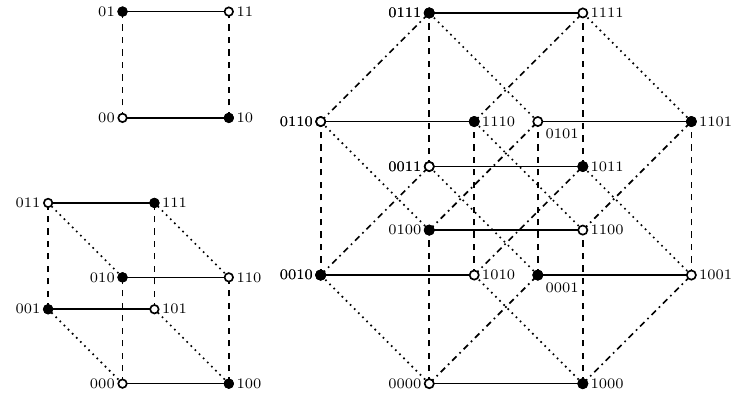}
\caption{The 2-, 3-, and 4-dimensional cubes.}
\label{Fig:Cubes}
\end{figure}

The {\bf $k$-skeleton} of $H_n$ is the union of all its $k$-faces, so the 2-skeleton is the
union of all the (square) faces.
The 1-skeleton is the {\bf hypercube graph}, denoted $Q_n$. Its
vertices are the $n$-digit binary strings, and an edge
connects any two vertices that differ in exactly one position.
Figure~\ref{Fig:Cubes} shows $Q_2$, $Q_3$ and $Q_4$. 

The hypercube graph $Q_n$ of any dimension is {\bf bipartite}, that is, its vertices can be partitioned into
two sets (say, black and white) such that each edge joins a black vertex to a white vertex.
(The black vertices are those binary strings with an odd number of 1's, while
the white vertices have an even number of 1's, as in Figure~\ref{Fig:Cubes}.)

Figure~\ref{Fig:Cubes} illustrates another nice feature of hypercubes:
$H_n$ can be edge-colored with colors $1,\ldots, n$, so that an edge
whose endpoints differ in the $k$th coordinate gets color $k$.
Note that each vertex is incident with exactly one edge of each color.

\section{Two-cell embeddings of graphs.}
\label{Section:Surfaces}

Surfaces are mathematical
objects that look locally like the plane, but which can
differ radically from the plane globally. The Earth's surface is an example:
it looks locally like a plane (at least on flat land) while globally it is not a plane at all, but a sphere.
More exactly, a {\bf surface} is a compact connected topological space in which each point has a neighborhood that is
homeomorphic to an open disk.
See introductory texts such as \cite{GOODMAN}
and \cite[Chapter 7]{CHARTRAND}
for development of the informal remarks we make here.

\begin{figure}[b]
\centering
\includegraphics{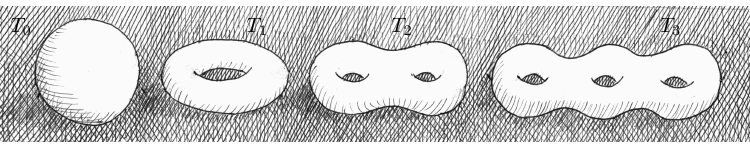}
\caption{Examples of Tori. The sphere $T_0$ (left) followed by $T_1$, $T_2$, and $T_3$.}
\label{Fig:Tori}
\end{figure}

The {\bf orientable} surfaces are the sphere, the torus,
the 2-holed torus, and, in general, surfaces with $g$ holes. (Nonorientable surfaces
are those that one can cut a M\"obius band out of, but they are not our focus here.)
The number of holes in an orientable surface is called its {\bf genus}.
We denote the (unique!) surface of genus~$g$ by $T_g$.
If $T$ is an arbitrary orientable surface,  its genus is denoted by $\gamma(T)$, so $\gamma(T_g)=g$.
The sphere has genus~0, and is thus denoted as $T_0$.
Figure~\ref{Fig:Tori} shows some examples.

As noted earlier,
the {\bf genus}  $\gamma(G)$ of a graph $G$ is the
smallest integer $g$ for which $G$ can be drawn on $T_g$
without crossed edges. Such a drawing of $G$ on $T_g$ is regarded as a continuous
injection $G\to T_g$, and is called an {\bf embedding} of $G$ in $T_g$.
 An embedding of $G$ in a surface of genus $\gamma(G)$
is called a {\bf genus embedding} of $G$.
Figure~\ref{Fig:K55} shows genus embeddings of $K_{3,3}$, $K_{4,4}$, and $K_{4,5}$.

An embedding $\varphi\colon G\to T$ divides $T$ into {\bf regions},
which are the connected components of $T \setminus \varphi(G)$.
A {\bf 2-cell embedding} is one in which each region is homeomorphic to
an open disk. 
The regions of a 2-cell embedding are called {\bf faces}.

Every genus embedding is a 2-cell embedding~\cite[Theorem 7.2]{CHARTRAND}. 
For example, the embedding of $K_{3,3}$ in Figure~\ref{Fig:K55} has three faces: two
squares and one octagon. The embedding of $K_{4,4}$ has eight faces, each a square.

{\bf Euler's formula}~\cite[Theorem 7.1]{CHARTRAND} implies
that if a 2-cell embedding of a connected graph in a closed orientable surface $T$ has $v$ vertices, $e$ edges, and $f$ faces, then
\begin{equation}
\gamma(T)=\frac{2-v+e-f}{2}.
\label{Eqn:Euler}
\end{equation}
Checking this on the embedding of $K_{4,4}$ in Figure~\ref{Fig:K55} 
we get $1=\frac{1}{2}(2-8+16-8)$. For $K_{4,5}$ in Figure~\ref{Fig:K55}, we get
$2=\frac{1}{2}(2-9+20-f)$, so there are $f=9$ faces.

If a 2-cell embedding of a connected graph with $v$ vertices and $e$ edges has $f$ faces
$F_1,\ldots, F_f$, and each $F_i$ is a $p_i$-gon,
then because each edge is on exactly two faces, we get $2e=p_1+\cdots +p_f$.
If the graph is bipartite, then it has no triangles, so $2e\geq 4f$.
Then equation~(\ref{Eqn:Euler}) yields a lemma~\cite[Corollary 7.6]{CHARTRAND}.

\begin{lemma}
If $G$ has $v$ vertices, $e$ edges, and is bipartite, then
$\displaystyle{\gamma(G)\geq\frac{1}{4}(4-2v+e)}$.
\label{Lemma}
\end{lemma}
\section{Genus embeddings of hypercubes.}
\label{Section:Genus}

We are ready for our theorem.

\begin{theorem}
The genus of the hypercube graph $Q_n$ is $\gamma(Q_n)=1+(n-4)2^{n-3}$.
\label{Theorem:Main}
\end{theorem}

\begin{figure}[H]
\centering
\includegraphics{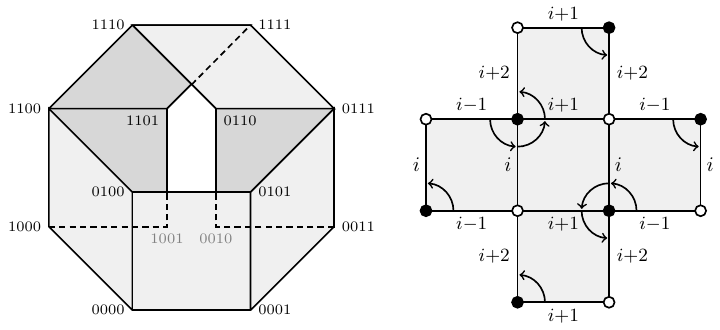}
\caption{Left:~A M\"obius strip in $H_4$. Right: Squares in $T$ are oriented by the
right-hand rule at black vertices. Squares bicolored $(i{-}1)i$,
$i(i{+}1)$, and $(i{+}1)(i{+}2)$ are shown.}
\label{Fig:Mobius}
\end{figure}

\begin{proof}
Color the edges of $H_n$ with colors $1,2,\ldots, n$, so that any edge joining vertices that
differ in the $k$th coordinate is given the color $k$. Any face of $H_n$ whose edges are
colored $k$ and $\ell$ is then
``bicolored'' by the pair $k\ell$. Assemble the collection $\mathcal{F}$ of all 
faces of $H_n$
having one of the bicolors $12,\; 23, \;34, \ldots,\, (n{-}1)n,\; n1$.
For any one of these $n$ bicolors, $H_n$ has $2^{n-2}$ faces of that particular bicolor,
so $|\mathcal{F}|=n2^{n-2}$.
Any edge $e$ of $Q_n$ belongs to exactly two faces in $\mathcal{F}$: if $e$ has
color $k$, then these two faces have bicolors $(k{-}1)k$ and $k(k{+}1)$
(addition modulo $n$). 
Thus $n$ faces of $\mathcal{F}$ are arranged cyclicly around each vertex,
in the manner described by Figure~\ref{Fig:Surface}.
It follows that the faces $\mathcal{F}$ form a surface
 $T$ in the 2-skeleton of $H_n$, with $Q_n$ embedded in it.
This surface has $n2^{n-2}$ square regions, and it is
connected because $Q_n$ is connected. 

Assume for the moment that this surface is orientable.
By the genus formula~(\ref{Eqn:Euler}),
\[\gamma(T)=\frac{2-v+e-f}{2}=\frac{2-2^n+n2^{n-1}-n2^{n-2}}{2}=1+(n-4)2^{n-3}.\]
Thus we have embedded $Q_n$ in a surface of genus $1+(n-4)2^{n-3}$.
Is this the lowest genus possible?
Lemma~\ref{Lemma} says yes: $\gamma(Q_n)\geq\frac{1}{4}(4-2\cdot 2^n+n2^{n-1})=\gamma(T)$.

 \begin{figure}[h]
 \centering
\includegraphics{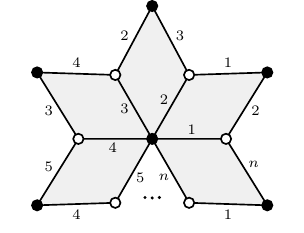}
 \caption{The squares in $\mathcal{F}$ that surround a vertex.}
 \label{Fig:Surface}
 \end{figure}
 
To finish the proof we must verify that $T$ is orientable. There is indeed something to prove
here, for when $n\,{\geq}\,4$, the 2-skeleton of $H_n$ contains M\"obius strips
(Figure~\ref{Fig:Mobius}, left). We must verify that none exist in~$T$.
To do this, note that each face of~$T$ has a local orientation given by the ``right-hand rule''
at either one of its black vertices: place your right hand at a black vertex with fingers
pointing from edge color $i$ to edge color $i{+}1$. Your thumb points in an ``up'' direction
for this square, and your other fingers indicate a counterclockwise orientation for the square.
(At white vertices, the thumb points ``down.'') Figure~\ref{Fig:Mobius} shows that this
orientation is preserved as we move from square to adjacent square on $T$.
It follows that $T$ is orientable.
\end{proof}

Let's carry out the construction of this proof for $Q_3$, $Q_4$, and $Q_5$.

\begin{figure}[H]
\centering
\includegraphics{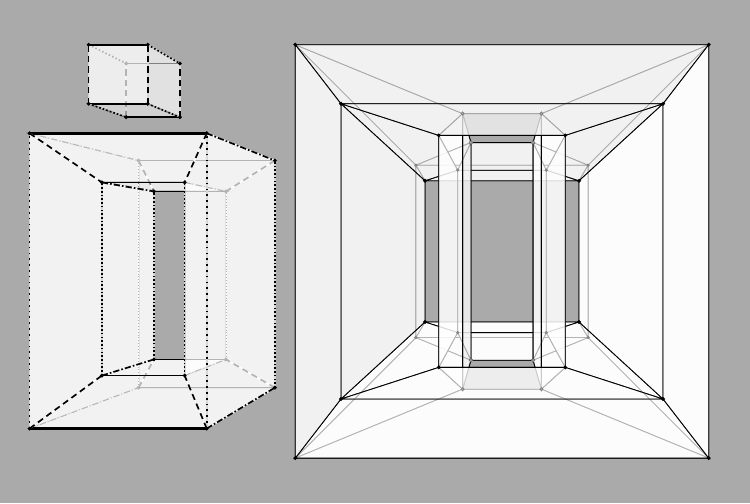}
\caption{Genus embeddings of $Q_3$, $Q_4$, and $Q_5$}
\label{Fig:Domus}
\end{figure}

First, $Q_3$. Color the edges of $Q_3$ solid, dashed, and dotted, as in Figure~\ref{Fig:Cubes}.
By our construction, we form a surface $T$ by including the solid/dashed faces, 
the dashed/dotted faces, and dotted/solid faces of $H_3$. These are in fact all the faces of
$H_3$, and we get the six faces shown in the upper left of Figure~\ref{Fig:Domus}.
They form the surface of the 3-cube, which is topologically equivalent to the sphere $T_0$.

Next, consider $Q_4$, with edges colored as in
Figure~\ref{Fig:Cubes}.
Our construction dictates that we form a surface $T$ by including the solid/dashed faces, 
the dashed/dotted faces, the dotted/dash-dotted faces, and the dash-dotted/solid
 faces of $H_4$. 
There are sixteen such faces. We see them on the bottom left of Figure~\ref{Fig:Domus}.
The resulting surface is $T_1$, the torus. This surface does not include the solid/dotted
 and dashed/dash-dotted faces
of $H_4$, but we clearly see their perimeters because their edges belong to the included faces.
In walking through the hole of the torus, one walks through all four solid/dotted faces and sees the
perimeters of all
four dashed/dash-dotted faces, which are inside the torus. 
The representation of $H_4$ in Figure~\ref{Fig:Cubes} has a long history.
According to Robbin~\cite{ROBBIN} such a perspective view of $H_4$ (which is now standard)
originated with V.~Schlegel (1843--1905).

Let's now apply our construction to get a genus embedding of $Q_5$. We
will not assign specific colors to the edges so as to avoid visual clutter in the image.
Say the five edge colors of $Q_5$ are 1, 2, 3, 4, and 5. We thus include the 5-cube faces
bicolored 12, 23, 34, 45, and 51 to obtain the embedding of $Q_5$ in $T_5$, shown
in Figure~\ref{Fig:Domus}.

This embedding does not include the 5-cube faces bicolored 13, 35, 52, 24, and 41. Thus the
embedding uses exactly half the faces of the 5-cube.  One can walk through this model and see all the edges and half the faces of $H_5$,
without any intersections. The other half of the faces of the 5-cube form a surface isometric to this one.  (See Section~\ref{Section:Parallel} below.) The missing faces are clearly visible because their
perimeters are edges of faces that do belong to the embedding. 

One nice feature of these embeddings is that they aid greatly in the visualization
of hypercubes.
One can, for example, easily pick out all ten facets of $Q_5$, and see how they fit together.
To highlight this, we offer a few exercises related to our model of $H_5$ in $T_5$
(Figure~\ref{Fig:Domus}).

\medskip
\noindent
{\bf Exercise 1:} Locate all 80 faces (squares) of $H_5$ in this model.

\smallskip
\noindent
{\bf Exercise 2:} Identify all ten facets of $H_5$. (The facets are 4-cubes.)

\smallskip
\noindent
{\bf Exercise 3:} The 5-cube has 40 3-faces (each one a 3-cube). Find them (or at least some of them). For each one, locate the two 4-cube facets that share it.

\smallskip
\noindent
{\bf Exercise 4:} For an arbitrary vertex, find five 4-cube facets that share this vertex.

\section{Parallel genus embeddings.}
\label{Section:Parallel}
In  \cite{RINGEL1} and ~\cite{BEINEKE}, the approach to cube embedding is {\it extrinsic};
they describe a recursive procedure that hooks together two lower-genus surfaces by a family of connecting ``tubes.''
In contrast, our {\it intrinsic} method has the interesting consequence that the entire 2-skeleton can be decomposed into copies of the genus surface.

Recall that a cycle $Z$ in a graph is a {\bf Hamiltonian cycle} if $Z$ contains
all  vertices of the graph. The {\bf complete graph} $K_n$ is the graph with
vertex set $\{1,2,\ldots,n\}$ and with an edge joining each pair of distinct vertices.
Any ordering  $i_1i_2\ldots i_n$ of $\{1,2,\ldots,n\}$ gives rise to a Hamiltonian cycle
$Z=i_1i_2\ldots i_n$ in $K_n$ whose $n$ edges are $\{i_ki_{k{+}1}\mid 1 \leq k \leq n\}$ (arithmetic modulo $n$).

Arguing as in the proof of Theorem~\ref{Theorem:Main},
if $Z = i_1i_2\ldots i_n$ is {\it any} Hamiltonian cycle in $K_n$, then the union
of the faces of $H_n$ that are bicolored $i_ki_{k+1}$ ($1 \leq k \leq n$), addition mod $n$, is
a surface $T(Z)$ that is a genus embedding for $Q_n\subseteq T(Z)$.  
In fact, if $\sigma$ is the permutation $k \mapsto i_k$, $1 \leq k \leq n$, then $\sigma$ induces a permutation of the standard unit vectors in $\mathbb{R}^n$ which carries $T$ isometrically onto $T(Z)$ sending vertices to vertices,
edges to edges, and faces to faces.
If $Z' =i_1'i_2'\ldots i_n'$ is another Hamiltonian cycle in $K_n$ and if $Z$ and $Z'$ share no edge, then $T(Z)\cap T(Z')=Q_n$.

We call a family of surfaces $T(Z_1), T(Z_2),\ldots,$ $T(Z_s)$ a {\bf parallel family} if every face (square) of
$H_n$ belongs to exactly one of the $T(Z_k)$.
A collection of Hamiltonian cycles $Z_1,Z_2,\ldots, Z_s$ of $K_n$ is
a {\bf Hamiltonian decomposition} of $K_n$ if each edge belongs to exactly one of the cycles. 

It has long been known that for odd $n$, $K_n$ has a Hamiltonian decomposition
into $(n{-}1)/2$ Hamiltonian cycles.  (This was a problem in recreational mathematics, asking whether one could find seating arrangements around a round table that gave each pair of people a unique side-by-side appearance. See~\cite{alspach}.)   We get a nice consequence for the hypercube.

\begin{theorem}
For odd $n \geq 3$, each Hamiltonian decomposition of $K_n$ gives rise to a parallel  family of $ (n{-}1)/2$ genus embeddings of $Q_n$. That is, the $2$-skeleton of $H_n$ can be decomposed into face-disjoint isometric copies of genus embeddings of $Q_n$.
\end{theorem}
\begin{proof}
Let $s = (n{-}1)/2$. Take a Hamiltonian decomposition $Z_1,\ldots,Z_s$ of $K_n$. 
The genus embeddings $T(Z_1),\ldots,T(Z_s)$
of $Q_n$ intersect pairwise in $Q_n$.
From the proof of Theorem \ref{Theorem:Main}, $T(Z_k)$ has $n 2^{n-2}$ faces, so the $s$  copies contain all $n2^{n-2}(n-1)/2={{n}\choose{2}} 2^{n-2}$ faces of $Q_n$.
\end{proof}


After completing this article, we found that Das~\cite{DAS} also
obtained Theorem 1 (but not Theorem 2) using a more general version of
our approach, but avoiding
the issue of nonorientability by induction on dimension. Das credits the idea
of using the 2-skeleton to Coxeter's constructions of certain skew
polyhedra~\cite{COXETER}. Indeed, two of Coxeter's skew polyhedra coincide with our genus surfaces for $Q_4$ and $Q_5$, though Coxeter does not refer to graph genus.

In fact, the above ideas can be generalized in other directions.  Any decomposition of $K_n$ into edge-disjoint cycles, possible for $n \geq 3$ odd by a theorem of Euler \cite[p. 64]{harary}, yields a parallel family of surfaces for the $2$-skeleton of $H_n$; see \cite{hk-ADAM}. But the parallel family given by Theorem 2 has all surfaces pairwise isometric.
Can this
polytopal perspective be applied to other graph genus questions? What 
about the {\it nonorientable} genus of the $n$-cube~\cite{JUNGERMAN}?  It seems that
many nice questions remain.

\begin{acknowledgment}{Acknowledgment.}
We thank the referees for remarks that led to significant improvements. RHH is supported by Simons Foundation Collaboration Grant for Mathematicians 523748.
\end{acknowledgment}

\begin{biog}
\item[Richard Hammack] 
is a professor of mathematics at Virginia Commonwealth University,
working mostly in graph theory and combinatorics. He is the author of
{\it Book of Proof}, an open proofs textbook, and coauthor of
{\it Handbook of Product Graphs} (with W.\ Imrich and S.\ Klav\v{z}ar).
He is a proud father of four daughters.
\begin{affil}
Department of Mathematics and Applied Mathematics, Virginia Commonwealth University,
Richmond,VA   23284-2014\\
rhammack@vcu.edu
\end{affil}

\end{biog}

\begin{biog}
\item[Paul Kainen] 
is adjunct professor of mathematics at Georgetown University. He is Director of the Lab for Visual Math at Georgetown and works primarily in graph theory, topology, and neural networks.
He published an editorial (Nov./Dec. 2016) in BioPhotonics, which proposed a strategy 
for laser-physics-based medicine, and he is coauthor of {\it The Four-Color Problem} (with T. L. Saaty).
\begin{affil}
Department of Mathematics and Statistics, Georgetown University,
Washington, DC 20057\\
kainen@georgetown.edu\\

\medskip
{\rm Citation of this article:\\

\smallskip

 R. H. Hammack \& P. C. Kainen, A new view of hypercube genus, {\it American Math. Monthly} {\bf 128} (4) (2021) 352--359.}
\end{affil}

\end{biog}

\vfill\eject


\begin{thebibliography}{99}

\bibitem{alspach}
Alspach, B. (2008). The wonderful Walecki construction. {\it Bull. Inst. Combin. Appl.} 52:  7--20.

\bibitem{BEINEKE}
Beineke, L. W., Harary F. (1965). The genus of the $n$-cube.
{\em Canad.~J.~Math}. 17:  494--496.

\bibitem{CHARTRAND} Chartrand, G., Lesniak L., Zhang, P. (2011).
{\it Graphs and Digraphs}, fifth edition. Boca Raton, FL: Chapman \& Hall/CRC Press. 

\bibitem{COXETER}
Coxeter, H. S. M.  (1937). Regular skew polyhedra in
three and four dimensions and their topological
anaolgues, {\it Proc. of the London Math. Soc.} 43(2): 33--62,


\bibitem{DAS}
Das, S. (2019). Genus of the hypercube graph and real moment-angle complexes,
{\it Topology and its Applications}. 258: 415--424.

\bibitem{GOODMAN} Goodman, S. E. (2005). {\it Beginning Topology}. Pacific Grove, CA: Thompson Brooks/Cole.

\bibitem{hk-ADAM}
Hammack, R. H., Kainen, P. C. (2020). Sphere decompositions of hypercubes. {\it Art of Discr. and Appl. Math.} Feb. 19, 2020.

\bibitem{harary}
Harary, F. (1969), {\it Graph Theory}. Reading, MA: Addison-Wesley.

\bibitem{HUNTER}
Hunter, R., Kainen, P. C.  (2007). Quadrilateral embedding of $G \times Q_s$.
{\it  Bull. Inst. Combin. Appl.} 52: 13--20.

\bibitem{JUNGERMAN} Jungerman, M.  (1978). The non-orientable genus of the $n$-cube.
{\it Pacific J. Math.} 76(2): 443--451.

\bibitem{pck-st-cr} Kainen, P. C. (1972). On the stable crossing number of cubes. {\it Proc. Amer. Math. Soc.} 36(1):  55--62.

\bibitem{k-white}
Kainen, P. C., White, A. T. (1978). On stable crossing numbers. {\it J.~Graph Theory}. 2(3): 181--187.

\bibitem{PISANSKI1}
Pisanski, T. (1992). Orientable quadrilateral embeddings of products of graphs. Algebraic graph theory (Leibnitz 1989), {\it Discrete Math.} 109(1-3): 203--205.

\bibitem{PISANSKI2}
Pisanski, T. (1980) Genus of Cartesian products of regular bipartite graphs {\it J. Graph Theory}.
4(1): 31--42.

\bibitem{RINGEL1} Ringel, G. (1955). \"Uber drei kombinatorische Probleme am $n$-dimensionalen
W\"urfel und W\"urfelgitter. {\it Abh. Math.~Sem.~Univ.~Hamburg}. 20: 10--19.

\bibitem{RINGEL2} Ringel, G. (1965). Das Gescshlecht des vollst\"andgen parren Graphen.
{\it Abh.~Math.~Sem.~Univ. Hamburg}. 28:  139--150.

\bibitem{RINGEL3} Ringel, G., Youngs, J. W. T. (1968). Solution of the Heawood map-coloring
problem. {\it Proc.~Nat.~Acad. Sci.~USA}. 60: 438--445.

\bibitem{ROBBIN} Robbin, T. (2006). {\it Shadows of Reality: The Fourth Dimension in Relativity,
Cubism and Modern Thought}. New Haven \& London: Yale University Press. 

\bibitem{WHITE} 
White, A. T. (1970). The genus of repeated cartesian products of bipartite graphs.
{\it Trans.~Amer.~Math.~Soc.} 151: 393--404; 

\bibitem{ZIEGLER}
Ziegler, G. M. (1994). {\it Lectures on Polytopes}. Series: Graduate Texts in Mathematics,
 152. New York: Springer Verlag.

\end{thebibliography}
\end{document}